\documentclass{scrartcl}
\usepackage{amsmath}
\usepackage{amsfonts}
\usepackage{amssymb}
\usepackage{dsfont}
\usepackage[T1]{fontenc}
\usepackage[latin1]{inputenc}
\usepackage{epsfig}
\usepackage{graphicx}

\usepackage[boxed, lined]{algorithm2e}
\usepackage{float}
\usepackage{comment}
\usepackage{mathtools}

\newcommand{\ba}{\mathbf{a}}
\newcommand{\bb}{\mathbf{b}}

\newcommand{\bw}{\mathbf{w}}
\newcommand{\bv}{\mathbf{v}}

\newcommand{\D}{{\mathcal{D}}}

\newcommand{\W}{{\mathcal{W}}}

\newcommand{\Nu}{{\mathcal{N}}}

\newcommand{\N}{\mathbb{N}}

\newcommand{\R}{\mathbb{R}}

\newcommand{\Rd}{\mathbb{R}^d}

\newcommand{\beq}{\begin{eqnarray*}}

\newcommand{\eeq}{\end{eqnarray*}}

\newcommand{\beqm}{\begin{eqnarray}}

\newcommand{\eeqm}{\end{eqnarray}}

\newtheorem{theorem}{Theorem}

\newtheorem{lemma}{Lemma}
\newtheorem{definition}{Definition}

\newcommand{\EXP}{{\mathbf E}}
\newcommand{\PROB}{{\mathbf P}}

\renewcommand{\P}{{\cal P}}
\renewcommand{\bf}{\normalfont \bfseries}
\renewcommand{\it}{\normalfont \itshape}

\allowdisplaybreaks
\begin{document}
\renewcommand{\thefootnote}{\fnsymbol{footnote}}
\newcommand{\F}{{\cal F}}
\newcommand{\Sp}{{\cal S}}
\newcommand{\G}{{\cal G}}
\newcommand{\HH}{{\cal H}}

\begin{center}

  {\LARGE \bf
    Analysis of the rate of convergence of an over-parametrized
    deep neural network estimate learned by gradient descent
  }
\footnote{
Running title: {\it Over-parametrized deep neural networks}}
\vspace{0.5cm}

Michael Kohler$^{1}$
and Adam Krzy\.zak$^{2,}$\footnote{Corresponding author. Tel:
  +1-514-848-2424 ext. 3007, Fax:+1-514-848-2830}\\

{\it $^1$
Fachbereich Mathematik, Technische Universit\"at Darmstadt,
Schlossgartenstr. 7, 64289 Darmstadt, Germany,
email: kohler@mathematik.tu-darmstadt.de}

{\it $^2$ Department of Computer Science and Software Engineering, Concordia University, 1455 De Maisonneuve Blvd. West, Montreal, Quebec, Canada H3G 1M8, email: krzyzak@cs.concordia.ca}

\end{center}
\vspace{0.5cm}

\begin{center}
September 1, 2022
\end{center}
\vspace{0.5cm}

\noindent
    {\bf Abstract}\\
Estimation of a regression function from independent and identically
distributed random variables is considered. The $L_2$ error with
integration
with respect to the design measure is used as an error criterion.
Over-parametrized deep neural network estimates are defined where
all the weights are learned by the gradient descent. It is shown that
the expected $L_2$ error of these estimates converges to zero with
the rate close to $n^{-1/(1+d)}$ in case that the
regression function is H\"older smooth with H\"older exponent
$p \in [1/2,1]$. In case of an interaction model where the
regression function is assumed to be a sum of H\"older smooth
functions where each of the functions depends   only on $d^*$ many
of $d$ components of the design variable, it is shown that these
estimates achieve the corresponding $d^*$-dimensional rate of convergence.

    \vspace*{0.2cm}

\noindent{\it AMS classification:} Primary 62G08; secondary 62G20.

\vspace*{0.2cm}

\noindent{\it Key words and phrases:}
neural networks,
nonparametric regression,
over-parametrization,
rate of convergence.

\section{Introduction}
\label{se1}
\subsection{Deep learning}
\label{se1sub1}
Deep learning, i.e., the fitting of deep neural networks to data,
has achieved tremendous success in various applications in the past ten years.
Deep neural networks are nowadays the most successful methods
in image classification (cf., e.g., Krizhevsky, Sutskever and Hinton  (2012)),
text classification
(cf., e.g., Kim (2014)),
machine translation
(cf., e.g., Wu et al. (2016))
or
mastering of games
(cf., e.g., Silver et al. (2017)).

Motivated by this huge success in applications there is also an
increasing
interest in theoretical properties of the estimates based on the deep neural
networks. Here in the past years various impressive results concerning
the rate of convergence of the least squares regression estimates based
on the deep neural networks have been derived (cf. , e.g.,
 Bauer and Kohler (2019),
Schmidt-Hieber (2020),
Kohler and Langer (2021), and the literature cited therein).
But these results ignore two important properties of the
deep neural network estimates applied in practice: Firstly,
the most successful estimates are usually over-parametrized
in the sense that the number of parameters of these estimates
is much larger than the sample size. And secondly, these estimates
are learned by the gradient descent applied to randomly initialized
weights of the neural network. This motivates the question:
If we define an over-parametrized
 deep neural network estimate by randomly initializing
its weights and by performing a suitable number of gradient
descent steps, does the resulting estimate then have nice
theoretical properties? The purpose of this article is to give results
that partially answer this question.

\subsection{Nonparametric regression}
\label{se1sub2}
We study deep neural networks in the context of nonparametric
regression. Here,
$(X,Y)$ is an
$\Rd \times \R$--valued random
vector $(X,Y)$ with $\EXP Y^2 < \infty$,
and $m(x)=\EXP\{Y|X=x\}$ is the corresponding
regression function $m:\Rd \rightarrow \R$.
Given a sample
of $(X,Y)$, i.e.,  a data set
\begin{equation}
  \label{se1eq1}
\D_n = \left\{
(X_1,Y_1), \ldots, (X_n,Y_n)
\right\},
\end{equation}
where
$(X,Y)$, $(X_1,Y_1)$, \ldots, $(X_n,Y_n)$ are i.i.d.,
the goal is to construct an estimator
\[
m_n(\cdot)=m_n(\cdot, \D_n):\Rd \rightarrow \R
\]
of the corresponding regression function
$m(x)=\EXP\{Y|X=x\}$ such that the so--called $L_2$ error
\[
\int |m_n(x)-m(x)|^2 {\PROB}_X (dx)
\]
is ``small'' (cf., e.g., Gy\"orfi et al. (2002)
for a systematic introduction to nonparametric regression and
a motivation for the $L_2$ error).

We are interested to investigate for given estimates $m_n$
how quickly the expected $L_2$ error
\begin{equation}
\label{se1eq2}
\EXP \int |m_n(x)-m(x)|^2 {\PROB}_X (dx)
\end{equation}
converges to zero.
It is well-known, that without regularity assumptions on the
smoothness
of $m$ it is not possible to derive nontrivial asymptotic bounds
on (\ref{se1eq2}) (cf.,
Theorem 7.2 and Problem 7.2 in
Devroye, Gy\"orfi and Lugosi (1996) and
Section 3 in Devroye and Wagner (1980)).
In order to formulate such regularity assumptions we will
use in this paper the notion of $(p,C)$--smoothness, which we
introduce next.

\begin{definition}
\label{se1de1}
  Let $p=q+s$ for some $q \in \N_0$ and $0< s \leq 1$.
A {\bf function} $m:\R^d \rightarrow \R$ is called
{\bf $(p,C)$-smooth}, if for every $\alpha=(\alpha_1, \dots, \alpha_d) \in
\N_0^d$
with $\sum_{j=1}^d \alpha_j = q$ the partial derivative
$\frac{
\partial^q m
}{
\partial x_1^{\alpha_1}
\dots
\partial x_d^{\alpha_d}
}$
exists and satisfies
\[
\left|
\frac{
\partial^q m
}{
\partial x_1^{\alpha_1}
\dots
\partial x_d^{\alpha_d}
}
(x)
-
\frac{
\partial^q m
}{
\partial x_1^{\alpha_1}
\dots
\partial x_d^{\alpha_d}
}
(z)
\right|
\leq
C
\cdot
\| x-z \|^s
\]
for all $x,z \in \R^d$, where $\Vert\cdot\Vert$ denotes the Euclidean norm.
\end{definition}
For $p \leq 1$ $(p,C)$--smoothness means that the function is
H\"older-smooth with exponent $p$ and H\"older-constant $C$.

\subsection{Main results}
\label{se1sub3}
Based on the previous work of the authors and their co-authors
we introduce in this article a technique to analyze the rate
of convergence of over-parametrized deep neural network
estimates learned by gradient descent. The key ingredients in this
theory are as follows: We control {\it generalization} ability
of the estimate by using a metric entropy bound of a class
of neural networks with bounded weights (cf., Lemma \ref{le4} below),
by imposing proper bounds on the weights during initialization,
and by choosing the number of gradient steps and the step size
properly. We analyze {\it optimization} of the empirical
$L_2$ risk during gradient descent by introducing a proper regularization,
and by optimization of the outer weights during the gradient
descent. Here we show at the same time that due to our restrictions
on the number of gradient descent steps and the step size the inner
weights do not change much (cf., Lemma \ref{le1} below). And we
control {\it approximation} by using the over-parametrization
and our special topology of the networks
to show that with high probability
a subset of the initial inner weights has nice properties and that this
in turn leads to good approximation properties of the corresponding
networks as soon as the outer weights are suitably chosen.

This new approach of analyzing over-parametrized deep neural
networks is illustrated by analyzing the rate of convergence
of an estimate introduced
in Section \ref{se2}
as follows: We choose topology of the network where the output of
the network
is defined as a linear combination of a huge number of fully connected
deep neural networks of constant widths and depth. We introduce special
initialization of the weights, where the output weights are zero and
all inner weights are generated with uniform distribution.
Then we perform a suitable large number of gradient descent steps with
a suitably small step size.
We show that the expected $L_2$ error of the truncated version of the
resulting estimate converges to zero with the rate of convergence close to
\[
n^{-\frac{1}{1+d}}
\]
in case that the regression function is $(p,C)$-smooth for some
$1/2 \leq p \leq 1$. Furthermore, we show that in case that
the regression function is the sum of H\"older smooth
functions where each of the functions depends   only of $d^*$
of $d$ components of $X$, our estimate achieves
the rate of convergence close to
\[
n^{-\frac{1}{1+d^*}}.
\]

\subsection{Discussion of related results}
\label{se1sub4}
In the last six years various results concerning the rate
of convergence of the least squares regression estimates
based on deep neural networks have been shown. One of the main
achievement in this area is derivation of good rates of
convergence for such estimates in case that the regression function
is a composition of functions where each of these functions
depends only on a few of its components. This was first shown in
Kohler and Krzy\.zak (2017) for $(p,C)$--smooth regression functions with
$p \leq 1$. Bauer and Kohler (2019) showed that such results also hold
in case $p>1$ provided the activation function of the network is
sufficiently smooth. The surprising fact that such results also
hold for the non-smooth ReLU activation function was shown in
Schmidt-Hieber (2020). Kohler and Langer (2021) proved that such result
can be also derived without imposing a sparsity constraint on the underlying
networks. Suzuki (2018) and Suzuki and Nitanda (2019) proved corresponding
results under weaker smoothness assumptions then in the above papers.
That the least squares estimates based on the deep neural networks
are able to adapt to some kind of local dimension of the regression function
was shown in Kohler, Krzy\.zak and Langer (2022)  and
Eckle and Schmidt-Hieber (2019). For further results on least squares
estimates based on deep neural networks we refer to
Imaizumi and Fukamizu (2019) and Langer (2021a)
and the literature cited therein.
Various approximation results for deep neural network can be found
in Lu et al. (2020), Yarotsky (2017),Yarotsky and Zhevnerchuk (2020)
and Langer (2021b).

For neural networks with one hidden layer Braun et al. (2021) analyzed
the gradient descent. There it was shown that in case of a proper
initialization estimates learned by the gradient descent can achieve a
rate of convergence of order $1/\sqrt{n}$ (up to a logarithmic factor)
in case that the Fourier transform of the regression functions decays
suitably fast. This decay of the Fourier transform is related to
the classical results of Barron (1993, 1994) for the least squares
estimate, where also a similar dimension-free
rate of convergence was proven in case
that the first moment of the Fourier transform of the regression function
is finite.
In case $d=1$ the above estimate of Braun et al. (2021) was analyzed
in Kohler and Krzy\.zak (2022) in an over-parametrized setting.
By controlling the complexity of the estimate via strong regularization
 it was possible in this article to show that
over-parametrization leads to an improved rate of convergence
in case $d=1$.

A review of various results on over-parametrized deep neural network
estimates learned by gradient descent can be found in
Bartlett, Montanari and Rakhlin (2021). These results usually
analyze the estimates in some asymptotically equivalent models
(like the mean field approach in
Mei, Montanari, and Nguyen (2018), Chizat and Bach (2018) or Nguyen and Pham (2020) or the neural tangent approach in Hanin and Nica (2019)). In contrast, we analyze
directly the expected $L_2$ error of the estimate in a standard
regression model.

It is well-known that gradient descent applied to deep neural network
can lead to estimates which minimizes the empirical $L_2$-risk
(cf., e.g., Allen-Zhu, Li and Song (2019), Kawaguchi and Huang (2019)
and the literature cited therein).
However, as was shown in Kohler and Krzy\.zak (2022) such estimates, in general, do not perform well on data, which is independent of the training data. In this article we avoid this problem by restricting
the number of gradient descent steps and the step sizes and by imposing
bounds on the absolute values of the initial weights.

Our approach is related to Drews and Kohler (2022), where the universal
consistence of over-parametrized deep neural network estimate
learned by gradient descent was shown. In fact, we generalize the results
there such that we are able to analyze the rate of convergence
of the estimates.Iin particular, this requires more precise analysis of the
approximation error, which we do by using a multiscale
approximation in Lemma \ref{le8} below.

Our bound on the covering number of over-parametrized
deep neural networks is based on the corresponding result
in Li, Gu and Ding (2021).

\subsection{Notation}
\label{se1sub5}
  The sets of natural numbers, real numbers and nonnegative real numbers
are denoted by $\N$, $\R$ and $\R_+$, respectively. For $z \in \R$, we denote
the smallest integer greater than or equal to $z$ by
$\lceil z \rceil$.
The Euclidean norm of $x \in \Rd$
is denoted by $\|x\|$.
For $f:\R^d \rightarrow \R$
\[
\|f\|_\infty = \sup_{x \in \R^d} |f(x)|
\]
is its supremum norm.
Let $\F$ be a set of functions $f:\Rd \rightarrow \R$,
let $x_1, \dots, x_n \in \Rd$, set $x_1^n=(x_1,\dots,x_n)$ and let
$p \geq 1$.
A finite collection $f_1, \dots, f_N:\Rd \rightarrow \R$
  is called an $L_p$ $\varepsilon$--cover of $\F$ on $x_1^n$
  if for any $f \in \F$ there exists  $i \in \{1, \dots, N\}$
  such that
  \[
  \left(
  \frac{1}{n} \sum_{k=1}^n |f(x_k)-f_i(x_k)|^p
  \right)^{1/p}< \varepsilon.
  \]
  The $L_p$ $\varepsilon$--covering number of $\F$ on $x_1^n$
  is the  size $N$ of the smallest $L_p$ $\varepsilon$--cover
  of $\F$ on $x_1^n$ and is denoted by $\Nu_p(\varepsilon,\F,x_1^n)$.

For $z \in \R$ and $\beta>0$ we define
$T_\beta z = \max\{-\beta, \min\{\beta,z\}\}$. If $f:\R^d \rightarrow
\R$
is a function and $\F$ is a set of such functions, then we set
$
(T_{\beta} f)(x)=
T_{\beta} \left( f(x) \right)$.

\subsection{Outline}
\label{se1sub6}
The over-parametrized deep neural network estimates considered
in this paper are introduced in Section \ref{se2}. The main results
are presented in Section \ref{se3}. Section \ref{se4} contains the proofs.

\section{Definition of the estimate}
\label{se2}

Throughout the paper we let
$\sigma(x)=1/(1+e^{-x})$
be the logistic squasher and
we define the topology of our neural networks as follows: We
let $K_n, L, r \in \N$ be parameters of our estimate and using
these parameters we
 set
\begin{equation}\label{se2eq1}
f_\bw(x) = \sum_{j=1}^{K_n} w_{1,1,j}^{(L)} \cdot f_{j,1}^{(L)}(x)
\end{equation}
for some $w_{1,1,1}^{(L)}, \dots, w_{1,1,K_n}^{(L)} \in \mathbb{R}$, where
$f_{j,1}^{(L)}=f_{\bw,j,1}^{(L)}$ are recursively defined by
\begin{equation}
  \label{se2eq2}
f_{k,i}^{(l)}(x) = \sigma\left(\sum_{j=1}^{r} w_{k,i,j}^{(l-1)}\cdot f_{k,j}^{(l-1)}(x) + w_{k,i,0}^{(l-1)} \right)
\end{equation}
for some $w_{k,i,0}^{(l-1)}, \dots, w_{k,i, r}^{(l-1)} \in \mathbb{R}$
$(l=2, \dots, L)$
and
\begin{equation}
  \label{se2eq3}
f_{k,i}^{(1)}(x) = \sigma \left(\sum_{j=1}^d w_{k,i,j}^{(0)}\cdot x^{(j)} + w_{k,i,0}^{(0)} \right)
\end{equation}
for some $w_{k,i,0}^{(0)}, \dots, w_{k,i,d}^{(0)} \in \mathbb{R}$.

This means that we consider neural networks which consist of $K_n$ fully
connected
neural networks of depth $L$ and width $r$ computed in parallel and compute
 a linear combination of the outputs of these $K_n$ neural
networks.
The weights in the $k$-th such network are denoted by
$(w_{k,i,j}^{(l)})_{i,j,l}$, where
$w_{k,i,j}^{(l)}$ is the weight between neuron $j$ in layer
$l$ and neuron $i$ in layer $l+1$.

We initialize the weights $\bw^{(0)}=((\bw^{(0)})_{k,i,j}^{(l))})_{k,i,j,l}$ as
follows: We set
\begin{equation}
\label{se2eq4}
(\bw^{(0)})_{1,1,k}^{(L)}=0
\quad (k=1, \dots, K_n),
\end{equation}
 we choose $(\bw^{(0)})_{k,i,j}^{(l)}$ uniformly distributed on
$[-c_1 \cdot (\log n)^2, c_1 \cdot (\log n)^2]$ if $l \in \{1, \dots, L-1\}$, and we
choose
$(\bw^{(0)})_{k,i,j}^{(0)}$ uniformly distributed on
$[-c_2 \cdot (\log n)^2 \cdot n^\tau, c_2 \cdot (\log n)^2\cdot n^\tau]$, where $\tau>0$ is a parameter of the estimate.
Here the random values are defined such that all components
of $\bw^{(0)}$ are independent.

After initialization of the weights we perform $t_n \in \N$ gradient
descent
steps each with a step size $\lambda_n>0$. Here we try to minimize
the regularized empirical $L_2$ risk
\begin{equation}
\label{se2eq5}
F_n(\bw)
= \frac{1}{n} \sum_{i=1}^n | Y_i - f_\bw (X_i)|^2
+ c_3 \cdot \sum_{k=1}^{K_n} |w_{1,1,k}^{(L)}|^2.
\end{equation}
 To do this we set
\begin{equation}
\label{se2eq6}
\bw^{(t)}
=
\bw^{(t-1)}
-
\lambda_n \cdot \nabla_{\bw} F_n(\bw^{(t-1)})
\quad
(t=1, \dots, t_n).
\end{equation}
Finally we define our estimate as a truncated version of the neural
network with weight vector $\bw^{(t_n)}$, i.e., we set
\begin{equation}
\label{se2eq7}
m_n(x)= T_{\beta_n} (f_{\bw^{(t_n)}}(x))
\end{equation}
 where $\beta_n = c_4 \cdot \log n$ and $T_{\beta} z
= \max\{ \min\{z, \beta\}, - \beta_n\}$ for $z \in \R$
and $\beta>0$.

\section{Main results}
\label{se3}

\subsection{A general theorem}
\label{se3sub1}

Our first result is a general theorem which we will apply
in the next two subsections in order to analyze the
rate of convergence of our over-parametrized deep neural
network estimate.

\begin{theorem}
  \label{th1}
Let $n \in \N$,
let $(X,Y)$, $(X_1,Y_n)$, \dots, $(X_n,Y_n)$
be independent and identically distributed $\Rd \times \R$--valued random variables such that $supp(X)$ is bounded and
\begin{equation}
\label{th1eq1}
\EXP\left\{
e^{c_5 \cdot Y^2}
\right\}
< \infty
\end{equation}
holds and that the corresponding regression function
$m(x)=\EXP\{Y|X=x\}$ is bounded.

Let $\sigma(x)=1/(1+e^{-x})$ be the logistic squasher, let
$K_n,L,r,t_n \in \N$, $\lambda_n, \tau>0$ and define the estimate
$m_n$ as in Section \ref{se2}.
 Let $\tilde{K}_n \in \{1, \dots, K_n\}$,
\[
w_{k,i,j}^{(l)} \in
[-c_1 \cdot (\log n)^2, c_1 \cdot (\log n)^2]
\quad (l=1, \dots, L, k=1, \dots \tilde{K}_n)
\]
and
\[
w_{k,i,j}^{(0)}
\in [-c_2 \cdot (\log n)^2\cdot n^\tau, c_2 \cdot (\log n)^2 \cdot n^\tau]
\quad (k=1, \dots, \tilde{K}_n).
\]
Assume
\begin{equation}
\label{th1eq*}
\left|
\sum_{k=1}^{\tilde{K}_n}
w_{1,1,k}^{(L)} \cdot f_{\bar{\bw},k,j}^{(L)}(x)
\right|
\leq \beta_n
\quad (x \in supp(X))
\end{equation}
for all $\bar{\bw}$ satisfying
$|\bar{w}_{i,j,k}^{(l)}-w_{i,j,k}^{(l)}| \leq \log n$
$(l=0, \dots, L-1)$.

Assume furthermore
\begin{equation}
\label{th1eq2}
\frac{K_n}{n^\kappa} \rightarrow 0
\quad (n \rightarrow \infty)
\end{equation}
for some $\kappa >0$ and
\begin{equation}
\label{th1eq3}
\frac{K_n}{\tilde{K}_n \cdot n^{ r \cdot (d+1) \cdot \tau+1}}
\rightarrow
\infty
\quad (n \rightarrow \infty),
\end{equation}
and that $t_n, \lambda_n$ are given by
 \begin{equation}
    \label{th1eq4}
    t_n = \lceil c_6 \cdot L_n \cdot \log n \rceil
\quad \mbox{and} \quad
\lambda_n= \frac{1}{L_n}
  \end{equation}
  for some $L_n>0$ which satisfies
  \begin{equation}
    \label{th1eq5}
    L_n \geq (\log n)^{10 \cdot L + 10} \cdot K_n^{3/2}.
  \end{equation}
Assume
\begin{equation}
  \label{th1eq6}
2 \cdot c_3 \cdot c_6 \geq 1 , \quad
c_4 \cdot c_5 \geq 1
\quad \mbox{and} \quad
4 \cdot c_4 \cdot c_6 \leq 1.
\end{equation}

Then we have for any $\epsilon >0$
\begin{eqnarray*}
&&
\EXP \int | m_n(x)-m(x)|^2 \PROB_X (dx)
\leq
c_7 \cdot
\Bigg(
\frac{ n^{\tau \cdot d + \epsilon}}{n} + \sum_{k=1}^{\tilde{K}_n} |w_{1,1,k}^{(L)}|^2
\\
&&
\hspace*{1cm}
+
\sup_{
(\bar{w}_{i,j,k}^{(l)})_{i,j,k,l} :
\atop
|\bar{w}_{i,j,k}^{(l)}-w_{i,j,k}^{(l)}| \leq \log n
\, (l=0, \dots, L-1)
}
\int
|
\sum_{k=1}^{\tilde{K}_n}
w_{1,1,k}^{(L)} \cdot f_{\bar{\bw},k,1}^{(L)}(x)-m(x)|^2 \PROB_X (dx)
\Bigg).
\end{eqnarray*}

 \end{theorem}

\noindent
    {\bf Remark 1.} The upper bound on the expected $L_2$ error
    of our neural network estimates corresponds to the usual
    bounds for least squares estimates (cf., e.g., Theorem 11.5
    in Gy\"orfi et al. (2002)). The term
    \[
\frac{ n^{\tau \cdot d + \epsilon}}{n}
\]
is used to bound the estimation error, which comes from the
fact that we minimize the empirical $L_2$ risk and not the $L_2$
risk during gradient descent. And
\[
\sup_{
(\bar{w}_{i,j,k}^{(l)})_{i,j,k,l} :
\atop
|\bar{w}_{i,j,k}^{(l)}-w_{i,j,k}^{(l)}| \leq \log n
\, (l=0, \dots, L-1)
}
\int
|
\sum_{k=1}^{\tilde{K}_n}
w_{1,1,k}^{(L)} \cdot f_{\bar{\bw},k,1}^{(L)}(x)-m(x)|^2 \PROB_X (dx)
\]
is used to bound the approximation error, which occurs since we
restrict our estimate to our class of neural networks. Observe that
here we can choose $\bw$ optimally in view of the upper
bound in Theorem \ref{th1}, so in fact our approximation error
also includes a minimum over all possible weights $\bw$.
The additional term
\[
\sum_{k=1}^{\tilde{K}_n} |w_{1,1,k}^{(L)}|^2
\]
is needed to bound an additional error due to gradient descent.

\noindent
    {\bf Remark 2.} Condition (\ref{th1eq6}) is e.g. satisfied, if we set
    \[
    c_4 = \frac{1}{c_5}, \quad
    c_6= \frac{c_5}{4} \quad \mbox{and} \quad
    c_3= \frac{1}{8 \cdot c_5}.
    \]

\subsection{Rate of convergence for $(p,C)$--smooth regression
  functions}
\label{se3sub2}

 The main challenge in deriving a rate
 of convergence from Theorem \ref{th1} above is to
 derive an approximation result for a smooth
 regression function and our neural networks such that the
 bounds on the weights of the networks are small. The first
 result in this respect is shown in our next theorem.

\begin{theorem}
  \label{th2}
Let $n \in \N$,
let $(X,Y)$, $(X_1,Y_n)$, \dots, $(X_n,Y_n)$
be independent and identically distributed $\Rd \times \R$
valued random variables which satisfy $supp(X) \subseteq [0,1]^d$ and
(\ref{th1eq1}).
Assume
that the corresponding regression function
$m(x)=\EXP\{Y|X=x\}$ is $(p,C)$-smooth for some
$1/2 \leq p \leq 1$ and some $C>0$.

Let $\sigma(x)=1/(1+e^{-x})$ be the logistic squasher, let
$L,r \in \N$ with $L  \geq 2$ and $r \geq 2d$, set
\[
K_n=n^{6d+r+2},
\]
\[
\tau= \frac{1}{1+d}
\]
and
\[
    t_n = \lceil c_6 \cdot L_n \cdot \log n \rceil
\quad \mbox{and} \quad
\lambda_n= \frac{1}{L_n}
\]
  for some $L_n>0$ which satisfies (\ref{th1eq5}).
  Define the estimate as in Section \ref{se2} and assume
  that (\ref{th1eq6}) holds.

Then we have for any $\epsilon >0$
\begin{eqnarray*}
&&
\EXP \int | m_n(x)-m(x)|^2 \PROB_X (dx)
\leq
c_8 \cdot n^{
- \frac{1}{1+d} + \epsilon
}.
\end{eqnarray*}

 \end{theorem}

\noindent
    {\bf Remark 3.} According to Stone (1982), the optimal
    minimax $L_2$ rate of convergence in case of $(p,C)$--smooth
    regression function is
    \[
n^{- \frac{2p}{2p+d}}.
    \]
    For $p=1/2$ our estimate achieves a rate of convergence, which
    is arbitrary close to this rate of convergence. For $p>1/2$ our
    derived rate of convergence is not optimal. We conjecture this is
    a consequence of our proof and not of property of the estimate, but it is
    an open problem to prove this.

\subsection{Rate of convergence in an interaction model}
\label{se3sub3}

In this subsection we assume that the regression function satisfies
\[
m(x)= \sum_{I \subseteq \{1, \dots, d\} \, : \, |I|=d^* } m_I (x_I),
\]
where $1 \leq d^* < d$, $m_I: \R^{d^*} \rightarrow \R$
$(I \subseteq \{1, \dots, d\}$, $ |I|=d^* )$ are $(p,C)$-smooth
functions and we use the notation
\[
x_I = (x^{(j_1)}, \dots, x^{(j_{d^*})})
\]
for $I=\{j_1, \dots, j_{d^*} \}$. Our aim is the modify the estimate
of Subsection \ref{se3sub2} such that it achieves in this case
the $d^*$-dimensional rate of convergence.

To achieve this, we define
\[
f_{\bw}(x) = \sum_{I \subseteq \{1, \dots, d\} \, : \, |I|=d^* } f_{\bw_I} (x_I)
\]
where $f_{\bw_I}$ is defined by (\ref{se2eq1})--(\ref{se2eq3}) with
$d$ replaced by $d^*$ and weight vector $\bw_I$, and
\[
\bw= \left( \bw_I \right)_{I \subseteq \{1, \dots, d\}, |I|=d^* }.
\]
We initialize the weights
$\bw^{(0)}=(((\bw^{(0)}_I)_{k,i,j}^{(l))})_{k,i,j,l})_{I \subseteq \{1, \dots, d\}, |I|=d^* }
$ as
follows: We set
\[
(\bw^{(0)}_I)_{1,1,k}^{(L)}=0
\quad (k=1, \dots, K_n,  I \subseteq \{1, \dots, d\},  |I|=d^* ),
\]
 we choose $(\bw^{(0)}_I)_{k,i,j}^{(l)}$ uniformly distributed on
$[-c_1 \cdot (\log n)^2, c_1 \cdot (\log n)^2]$ if $l \in \{1, \dots, L-1\}$, and we
choose
$(\bw^{(0)}_I)_{k,i,j}^{(0)}$ uniformly distributed on
$[-c_2 \cdot (\log n)^2 \cdot n^\tau, c_2 \cdot (\log n)^2\cdot
n^\tau]$, where $\tau>0$ is a parameter of the estimate defined in
Theorem \ref{th3} below
$(I \subseteq \{1, \dots, d\}, |I|=d^*)$.
Here the random values are defined such that all components
of $\bw^{(0)}$ are independent.

After initialization of the weights we perform $t_n \in \N$ gradient
descent
steps each with a step size $\lambda_n>0$. Here we try to minimize
the regularized empirical $L_2$ risk
\[
F_n(\bw)
= \frac{1}{n} \sum_{i=1}^n | Y_i - f_\bw (X_i)|^2
+ c_3 \cdot
\sum_{I \subseteq \{1, \dots, d\} \, : \, |I|=d^* }
\sum_{k=1}^{K_n} |(\bw_I)_{1,1,k}^{(L)}|^2.
\]
 To do this we set
\[
\bw^{(t)}
=
\bw^{(t-1)}
-
\lambda_n \cdot \nabla_{\bw} F_n(\bw^{(t-1)})
\quad
(t=1, \dots, t_n).
\]
Finally we define our estimate as a truncated version of the neural
network with weight vector $\bw^{(t_n)}$, i.e., we set
\[
m_n(x)= T_{\beta_n} (f_{\bw^{(t_n)}}(x))
\]
 where $\beta_n = c_4 \cdot \log n$

\begin{theorem}
  \label{th3}
Let
$d \in \N$, $d^* \in \{1, \dots, d\}$, $1/2 \leq p \leq 1$, $C>0$,
let $n \in \N$,
let $(X,Y)$, $(X_1,Y_n)$, \dots, $(X_n,Y_n)$
be independent and identically distributed $\Rd \times \R$
valued random variables such that $supp(X) \subseteq [0,1]^d$  and
(\ref{th1eq1})
holds. Assume that the corresponding regression function
$m(x)=\EXP\{Y|X=x\}$ satisfies
\[
m(x)= \sum_{I \subseteq \{1, \dots, d\} \, : \, |I|=d^* } m_I (x_I)
\quad (x \in [0,1]^d)
\]
for some $(p,C)$--smooth functions
$m_I:\R^{d^*} \rightarrow
\R$
$(I \subseteq \{1, \dots, d\}, |I|=d^*)$.

Let $\sigma(x)=1/(1+e^{-x})$ be the logistic squasher, let
$L,r \in \N$ with $L  \geq 2$ and $r \geq 2d^*$, set
\[
K_n=n^{6d^*+r+2},
\]
\[
\tau= \frac{1}{1+d^*}
\]
and
\[
    t_n = \lceil c_6 \cdot L_n \cdot \log n \rceil
\quad \mbox{and} \quad
\lambda_n= \frac{1}{L_n}
\]
  for some $L_n>0$ which satisfies (\ref{th1eq5}).
Define the estimate as above.

Then we have for any $\epsilon >0$
\begin{eqnarray*}
&&
\EXP \int | m_n(x)-m(x)|^2 \PROB_X (dx)
\leq
c_9 \cdot n^{
- \frac{1}{1+d^*} + \epsilon
}
.
\end{eqnarray*}

 \end{theorem}

\noindent
    {\bf Remark 4.}
    The rate of convergence derived in Theorem \ref{th3} does not depend
    on $d$, hence under the above assumption on the regression function
    our estimate is able to circumvent the curse of dimensionality.
    That this is possible is well-known (cf., Stone (1994) and the literature
    cited therein), however our result is the first result
    which shows that
    this is also possible for (over-parametrized) neural network estimates
    learned by the gradient descent.

    \section{Proofs}
\label{se4}

    \subsection{Auxiliary results for the proof of Theorem \ref{th1}}
    \label{se4sub1}
    In this section we present five auxiliary results which we will use
    in the proof of Theorem \ref{th1}. Our first auxiliary result
    will play key role in the analysis of gradient descent.

\begin{lemma}
  \label{le1}
  Let
  $F:\R^K \rightarrow \R_+$
  be a nonnegative differentiable function.
  Let
  $t \in \N$, $L>0$, $\ba_0 \in \R^K$ and set
  \[
  \lambda=
\frac{1}{L}
\]
and
\[
\ba_{k+1}=\ba_k - \lambda \cdot (\nabla_{\ba} F)(\ba_k)
\quad
(k \in \{0,1, \dots, t-1\}).
\]
Assume
\begin{equation}
  \label{le1eq1}
  \left\|
 (\nabla_{\ba} F)(\ba)
  \right\|
  \leq
  \sqrt{
2 \cdot t \cdot L \cdot \max\{ F(\ba_0),1 \}
    }
\end{equation}
for all $\ba \in \R^K$ with
$\| \ba - \ba_0\| \leq \sqrt{2 \cdot t \cdot \max\{ F(\ba_0),1 \} / L}$,
and
\begin{equation}
  \label{le1eq2}
\left\|
(\nabla_{\ba} F)(\ba)
-
(\nabla_{\ba} F)(\bb)
  \right\|
  \leq
  L \cdot \|\ba - \bb \|
\end{equation}
for all $\ba, \bb \in \R^K$ satisfying
\begin{equation}
  \label{le1eq3}
  \| \ba - \ba_0\| \leq \sqrt{8 \cdot \frac{t}{L} \cdot \max\{ F(\ba_0),1 \}}
  \quad \mbox{and} \quad
  \| \bb - \ba_0\| \leq \sqrt{8 \cdot \frac{t}{L} \cdot \max\{ F(\ba_0),1 \}}.
\end{equation}
Then we have
\[
\|\ba_k-\ba_0\| \leq
\sqrt{
2 \cdot \frac{k}{L} \cdot (F(\ba_0)-F(\ba_k))
}
\quad
 \mbox{for all }
 k \in \{1, \dots,t\},
\]
\[
\sum_{k=0}^{s-1}
\| \ba_{k+1}-\ba_k \|^2
\leq
\frac{2}{L}
 \cdot (F(\ba_0)-F(\ba_s))
\quad
 \mbox{for all }
 s \in \{1, \dots,t\}
 \]
 and
 \[
 F(\ba_k)
 \leq
 F(\ba_{k-1})
                -
                \frac{1}{2 L} \cdot
                \| \nabla_{\ba}  F(\ba_{k-1}) \|^2
                \quad
 \mbox{for all }
 k \in \{1, \dots,t\}.
\]
\end{lemma}

\noindent
    {\bf Proof.} The result follows from Lemma 2 in Braun et al. (2021)
    and its proof.
    \hfill $\Box$

    Our next auxiliary result will help us to show that assumption
    (\ref{le1eq1}) is satisfied in the proof of Theorem \ref{th1}.

\begin{lemma}\label{le2}
  Let $\sigma: \R \rightarrow \R$ be bounded and differentiable, and assume that
its derivative is bounded.
Let $\alpha_n \geq 1$,
$t_n \geq L_n$,
$\gamma_n^* \geq 1$, $B_n \geq 1$, $r \geq 2d$,
\begin{equation}
	\label{le2eq1}
	|w_{1,1,k}^{(L)}| \leq \gamma_n^* \quad (k=1, \dots, K_n),
\end{equation}
\begin{equation}
	\label{le2eq2}
	|w_{k,i,j}^{(l)}| \leq B_n
	\quad
	\mbox{for } l=1, \dots, L-1
\end{equation}
and
\begin{equation}
	\label{le2eq3}
	\|\bw-\bv\|_\infty^2 \leq \frac{2t_n}{L_n} \cdot \max\{ F_n(\bv),1 \}.
\end{equation}
Then we have
\[
\| (\nabla_\bw F_n)(\bw) \|
\leq
c_{10} \cdot K_n^{3/2} \cdot B_n^{2L} \cdot (\gamma_n^*)^2 \cdot \alpha_n^{2} \cdot \sqrt{\frac{t_n}{L_n} \cdot \max\{F_n(\bv),1\}}.
\]
\end{lemma}

\noindent
    {\bf Proof.}
    See Lemma 2 in Drews and Kohler (2022).
    \hfill $\Box$

    Our third auxiliary result will help us to show that assumption
    (\ref{le1eq2}) is satisfied in the proof of Theorem \ref{th1}.

    \begin{lemma}
      \label{le3}
     Let $\sigma: \R \rightarrow \R$ be bounded and differentiable, and assume that its derivative
     is
     Lipschitz continuous and bounded.
     Let $\alpha_n \geq 1$,
     $t_n \geq L_n$,
     $\gamma_n^* \geq 1$, $B_n \geq 1$, $r \geq 2d$ and assume
     \begin{equation}
     	\label{le3eq1}
     	|\max\{ (\bw_1)_{1,1,k}^{(L)}, (\bw_2)_{1,1,k}^{(L)}\}| \leq \gamma_n^* \quad (k=1,
     	\dots, K_n),
     \end{equation}

     \begin{equation}
     	\label{le3eq2}
     	|\max\{(\bw_1)_{k,i,j}^{(l)},(\bw_2)_{k,i,j}^{(l)}\}| \leq B_n
     	\quad
     	\mbox{for } l=1, \dots, L-1
     \end{equation}
     and
     \begin{equation}
     	\label{le3eq3}
     	\|\bw_2-\bv\|^2 \leq 8 \cdot \frac{t_n}{L_n} \cdot \max\{ F_n(\bv),1 \}.
     \end{equation}
     Then we have
     \begin{eqnarray*}
     	&&
     	\| (\nabla_\bw F_n)(\bw_1) - (\nabla_\bw F_n)(\bw_2) \| \\
     	&&
     	\leq
     	c_{11} \cdot \max \{\sqrt{F_n(\bv)},1\} \cdot (\gamma_n^*)^{2} \cdot B_n^{3L-1} \cdot \alpha_n^{3} \cdot K_n^{3/2} \cdot \sqrt{\frac{t_n}{L_n}} \cdot \|\bw_1-\bw_2\|.
     \end{eqnarray*}
\end{lemma}

\noindent
    {\bf Proof.}     See Lemma 3 in Drews and Kohler (2022).
    \hfill $\Box$

    Our fourth auxiliary result uses a metric entropy bound in order
    to control the complexity of a set of over-parametrized deep
    neural networks.

\begin{lemma}
  \label{le4}
  Let $\alpha \geq 1$, $\beta>0$ and let $A,B,C \geq 1$.
  Let $\sigma:\R \rightarrow \R$ be $k$-times differentiable
  such that all derivatives up to order $k$ are bounded on $\R$.
  Let $\F$
  be the set of all functions $f_{\bw}$ defined by
  (\ref{se2eq1})--(\ref{se2eq3}) where the weight vector $\bw$
  satsifies
  \begin{equation}
    \label{le4eq1}
    \sum_{j=1}^{K_n} |w_{1,1,j}^{(L)}| \leq C,
    \end{equation}
  \begin{equation}
    \label{le4eq2}
    |w_{k,i,j}^{(l)}| \leq B \quad (k \in \{1, \dots, K_n\},
    i,j \in \{1, \dots, r\}, l \in \{1, \dots, L-1\})
    \end{equation}
and
  \begin{equation}
    \label{le4eq3}
    |w_{k,i,j}^{(0)}| \leq A \quad (k \in \{1, \dots, K_n\},
    i \in \{1, \dots, r\}, j \in \{1, \dots,d\}).
  \end{equation}
  Then we have for any $1 \leq p < \infty$, $0 < \epsilon < \beta$ and
  $x_1^n \in [-\alpha,\alpha]^d$
  \begin{eqnarray*}
    &&
  \Nu_p \left(
\epsilon, \{ T_\beta f  \, : \, f \in \F \}, x_1^n
\right)
\\
&&
\leq
\left(
c_{12} \cdot \frac{\beta^p }{\epsilon^p}
\right)^{
c_{13} \cdot \alpha^{d} \cdot A^{d} \cdot B^{(L-1) \cdot d} \left(\frac{C}{\epsilon}\right)^{d/k} + c_{14}
  }.\\
  \end{eqnarray*}

  \end{lemma}

\noindent
    {\bf Proof.}
Standard application of the chain rule for derivation together with
the
above bounds on the weight vector $\bw$
shows that we have for any
$x \in \Rd$ and any $s_1, \dots, s_k \in \{1, \dots, d\}$
\[
\left|
\frac{
\partial^k f_\bw
}{
\partial x^{(s_1)} \dots \partial x^{(s_k)}
}
(x)
\right|
\leq
c_{15} \cdot C \cdot B^{(L-1) \cdot k} \cdot A^k=:c
\]
(cf., e.g., proof of
 Lemma 4 in Drews and Kohler (2022)).

Let $\G \circ \Pi$ be the set of all piecewise polynomials
of total degree less than $k$ with respect to a partition $\Pi$
of $[-\alpha,\alpha]^d$ into cubes of sidelength
\[
\left(
c_{16} \cdot
\frac{\epsilon}{c}
\right)^{1/k},
\]
where $c_{16}=c_{16}(d,k)$ is a suitable small constant greater zero.
Then a standard bound on the remainder of a multivariate Taylor
polynomial
shows that for each $f_{\bw}$ we can find $g \in \G \circ \Pi$ such
that
\[
|f_{\bw}(x) - g(x)| \leq \frac{1}{2}
\]
holds for all $x \in [-\alpha,\alpha]^d$, which implies
\[
  \Nu_p \left(
\epsilon, \{ T_\beta f  \, : \, f \in \F \}, x_1^n
\right)
\leq
  \Nu_p \left(
\frac{\epsilon}{2}, \{ T_\beta g  \, : \, g \in \G \}, x_1^n
\right).
\]
$\G$ is a linear vector space of dimension less than or equal to
\[
c_{17} \cdot \alpha^d \cdot \left(
\frac{c}{\epsilon}
\right)^{d/k},
\]
from which we get the assertion by an application of Theorems 9.4 and 9.5
in Gy\"orfi et al. (2002).
    \hfill $\Box$

In order to be able to formulate our next auxiliary result we
need the following notation:
Let $(x_1,y_1), \dots, (x_n,y_n) \in \Rd \times \R$, let $K \in \N$,
let $B_1,\dots,B_K:\Rd \rightarrow \R$ and let $c_3>0$. In the next lemma
we consider the problem to minimize
\begin{equation}
  \label{se5eq1}
  F(\ba) =
  \frac{1}{n} \sum_{i=1}^n
  |\sum_{k=1}^K a_k \cdot B_k(x_i)-y_i|^2
  +
    c_3 \cdot  \sum_{k=1}^{K_n} a_k^2 ,
  \end{equation}
where $\ba=(a_1,\dots,a_K)^T$,
by gradient descent. To do this, we choose $\ba^{(0)} \in \R^K$
and set
\begin{equation}
  \label{se5eq2}
  \ba^{(t+1)} = \ba^{(t)}
  - \lambda_n \cdot (\nabla_\ba F)(\ba^{(t)})
    \end{equation}
for some properly chosen $\lambda_n>0.$

        \begin{lemma}
          \label{le6}
          Let $F$ be defined by (\ref{se5eq1}) and choose $\ba_{opt}$
          such that
          \[
F(\ba_{opt})=\min_{\ba \in \R^{K}} F(\ba).
          \]
          Then for any
          $\ba \in \R^{K}$
          we have
                    \[
                            \|(\nabla_\ba F)(\ba)\|^2
                            \geq 4 \cdot c_3  \cdot (F(\ba)-F(\ba_{opt})).
                            \]
          \end{lemma}

        \noindent
            {\bf Proof.}
 See Lemma 8  in Drews and Kohler (2022).
    \hfill $\Box$

        \subsection{Proof of Theorem \ref{th1}}
    \label{se4sub2}
In the proof we combine ideas from the proof of Theorem 1 in Drews and
Kohler (2022)  with ideas from the proof of Lemma 1 in Bauer and
Kohler (2019).

W.l.o.g. we assume
throughout the proof that $n$ is sufficiently large and that
$\|m\|_\infty \leq \beta_n$ holds.
   Let $A_n$ be the event that firstly the weight vector $\bw^{(0)}$
            satisfies
            \[
            | (\bw^{(0)})_{j_s,k,i}^{(l)}-\bw_{j_s,k,i}^{(l)}| \leq \log n
            \quad \mbox{for all } l \in \{0, \dots, L-1\},
            s \in \{1, \dots, \tilde{K}_n \}
            \]
            for some pairwise distinct $j_1, \dots, j_{\tilde{K}_n}
            \in \{1, \dots, K_n\}$
and such that secondly
\[
\max_{i=1, \dots, n} |Y_i| \leq \sqrt{\beta_n}
\]
holds.

Define the weight vectors
$(\bw^*)^{(t)}$ by
\[
((\bw^*)^{(t)})_{k,i,j}^{(l)} = (\bw^{(t)})_{k,i,j}^{(l)} \quad
\mbox{for all } l=0,\dots, L-1
\]
and
\[
((\bw^*)^{(t)})_{1,1,j_k}^{(L)} = w_{1,1,k}^{(L)} \quad \mbox{for all } k=1,\dots, \tilde{K}_n
\]
and
\[
((\bw^*)^{(t)})_{1,1,k}^{(L)} = 0 \quad \mbox{for all } k \notin \{j_1,\dots,
j_{\tilde{K}_n}\}.
\]

We
decompose the  $L_2$ error of $m_n$ in a sum of several terms.
Set
\[
m_{\beta_n}(x)=\EXP\{ T_{\beta_n} Y | X=x \}.
\]
We have
\begin{eqnarray*}
&&
\int | m_n(x)-m(x)|^2 \PROB_X (dx)
\\
&&
=
\left(
\EXP \left\{ |m_n(X)-Y|^2 | \D_n \right\}
-
\EXP \{ |m(X)-Y|^2\}
\right)
\cdot 1_{A_n}
+
\int | m_n(x)-m(x)|^2 \PROB_X (dx)
\cdot 1_{A_n^c}
\\
&&
=
\Big[
\EXP \left\{ |m_n(X)-Y|^2 | \D_n \right\}
-
\EXP \{ |m(X)-Y|^2\}
\\
&&
\hspace*{2cm}
- \left(
\EXP \left\{ |m_n(X)-T_{\beta_n} Y|^2 | \D_n \right\}
-
\EXP \{ |m_{\beta_n}(X)- T_{\beta_n} Y|^2\}
\right)
\Big] \cdot 1_{A_n}
\\
&&
\quad +
\Big[
\EXP \left\{ |m_n(X)-T_{\beta_n} Y|^2| \D_n \right\}
-
\EXP \{ |m_{\beta_n}(X)- T_{\beta_n} Y|^2\}
\\
&&
\hspace*{2cm}
-
2 \cdot \frac{1}{n} \sum_{i=1}^n
\left(
|m_n(X_i)-T_{\beta_n} Y_i|^2
-
|m_{\beta_n}(X_i)- T_{\beta_n} Y_i|^2
\right)
\Big] \cdot 1_{A_n}
\\
&&
\quad
+\Big[
2 \cdot \frac{1}{n} \sum_{i=1}^n
|m_n(X_i)-T_{\beta_n} Y_i|^2
-
2 \cdot \frac{1}{n} \sum_{i=1}^n
|m_{\beta_n}(X_i)- T_{\beta_n} Y_i|^2
\\
&&
\hspace*{2cm}
- \left(
2 \cdot \frac{1}{n} \sum_{i=1}^n
|m_n(X_i)-Y_i|^2
-
2 \cdot \frac{1}{n} \sum_{i=1}^n
|m(X_i)- Y_i|^2
\right)
\Big] \cdot 1_{A_n}
\\
&&
\quad
+
\Big[
2 \cdot \frac{1}{n} \sum_{i=1}^n
|m_n(X_i)-Y_i|^2
-
2 \cdot \frac{1}{n} \sum_{i=1}^n
|m(X_i)- Y_i|^2
\Big] \cdot 1_{A_n}
\\
&&
\quad
+
\int | m_n(x)-m(x)|^2 \PROB_X (dx)
\cdot 1_{A_n^c}
\\
&&
=: \sum_{j=1}^5 T_{j,n}.
\end{eqnarray*}
In the reminder of the proof we bound
\[
\EXP T_{j,n}
\]
for $j \in \{1, \dots, 5\}$.

In the {\it first step of the proof} we show
\[
\EXP T_{1,n} \leq c_{18} \cdot \frac{\log n}{n}.
\]
This follows as in the proof of Lemma 1 in Bauer and Kohler (2019).

In the {\it second step of the proof} we show
\[
\EXP T_{3,n} \leq c_{19} \cdot \frac{\log n}{n}.
\]
Again this follows from the proof of Lemma 1 in Bauer and Kohler (2019).

In the {\it third step of the proof} we show
\[
\EXP T_{5,n} \leq c_{20} \cdot \frac{(\log n)^2}{n}.
\]
The definition of $m_n$ implies $\int |m_n(x)-m(x)|^2 \PROB_X (dx) \leq
4 \cdot c_4^2 \cdot (\log
n)^2$, hence it suffices to show
\begin{equation}
\label{pth1eq1}
\PROB(A_n^c) \leq \frac{c_{21}}{n}.
\end{equation}
To do this, we consider sequential choice of the weights of the
$K_n$ fully connected neural networks. Probability that the weights in
the first of these networks differ in all components at most by $\log n$
from $w_{1,i,j}^{(l)}$ $(l=0, \dots, L-1)$ is
for large $n$ bounded below by
\begin{eqnarray*}
 \left( \frac{\log n}{2 \cdot c_1 \cdot (\log n)^2}
\right)^{r \cdot (r+1) \cdot (L-1)}
\cdot
\left(
\frac{\log n}{2 \cdot c_2 \cdot n^\tau}
\right)^{r \cdot (d+1)}
&\geq&
 n^{- r \cdot (d+1) \cdot \tau-0.5}.
\end{eqnarray*}
Hence probability that none of the first $n^{r \cdot (d+1) \cdot
  \tau +1}$ neural networks satisfies this condition is for large $n$
 bounded above by
\begin{eqnarray*}
(1 -  n^{- r \cdot (d+1) \cdot \tau-0.5}) ^{n^{r \cdot (d+1) \cdot
  \tau +1}}
&\leq&
\left(\exp \left(
-  n^{- r \cdot (d+1) \cdot \tau-0.5}
\right)
\right) ^{n^{r \cdot (d+1) \cdot
  \tau +1}}
\\
&=&
\exp( -  n^{0.5}).
\end{eqnarray*}
Since we have $K_n \geq n^{r \cdot (d+1) \cdot
  \tau +1} \cdot \tilde{K}_n$ for $n$ large we can successively
use the same construction for all of $\tilde{K}_n$ weights and we can conclude:
Probability that there exists $k \in \{1, \dots, \tilde{K}_n\}$
such that
none of the $K_n$ weight vectors of the fully
connected neural network differs by at most $\log n$ from
$(w_{i,j,k}^{(l)})_{i,j,l}$ is for large $n$ bounded from above by
\begin{eqnarray*}
&&
\tilde{K}_n \cdot \exp( -  n^{0.5})
\leq  n^\kappa \cdot \exp( - \cdot n^{0.5})
\leq \frac{c_{22}}{n}.
\end{eqnarray*}
This implies for large $n$
\begin{eqnarray*}
\PROB(A_n^c)
&\leq&
\frac{c_{22}}{n}
+
\PROB\{ \max_{i=1, \dots, n} |Y_i| > \sqrt{\beta_n}
\}
 \leq
\frac{c_{22}}{n}
+
n \cdot\PROB\{  |Y| > \sqrt{\beta_n}
\}
\\
& \leq &
\frac{c_{22}}{n}
+
n \cdot
\frac{\EXP\{ \exp(c_5 \cdot Y^2)}{ \exp( c_5 \cdot \beta_n)}
\leq
\frac{c_{19}}{n},
\end{eqnarray*}
where the last inequlity holds because of (\ref{th1eq1})  and
$c_4 \cdot c_5 \geq 1$.

Let $\epsilon >0$ be arbitrary.
In the {\it fourth step of the proof} we show
\[
\EXP T_{2,n} \leq
c_{23} \cdot
\frac{ n^{\tau \cdot d + \epsilon}}{n}
.
\]

Let $\W_n$ be the set of all weight vectors
$(w_{i,j,k}^{(l)})_{i,j,k,l}$ which satisfy
\[
| w_{1,1,k}^{(L)}| \leq (c_1+1) \cdot (\log n)^2 \quad (k=1, \dots, K_n),
\]
\[
|w_{i,j,k}^{(l)}| \leq (c_1+1) \cdot (\log n)^2 \quad (l=1, \dots, L-1)
\]
and
\[
|w_{i,j,k}^{(0)}| \leq (c_2+1) \cdot (\log n)^2 \cdot n^\tau.
\]
By  Lemma \ref{le1}, Lemma \ref{le2} and Lemma  \ref{le3} we can
conclude that on $A_n$ we have
\begin{equation}
\label{pth1eq2}
\| \bw^{(t)}-\bw^{(0)}\| \leq \log n \quad (t=1,
\dots, t_n).
\end{equation}
This follows from the fact that on $A_n$ we have
\[
F_n(\bw^{(0)})
=
\frac{1}{n} \sum_{i=1}^n Y_i^2 \leq \beta_n
\]
and that
\[
\frac{2 \cdot t_n}{L_n} \cdot \beta_n
\leq
4 \cdot c_4 \cdot c_6 \cdot (\log n)^2
\leq (\log n)^2.
\]
Together with the initial choice of $\bw^{(0)}$ this implies that on
$A_n$
we have
\[
\bw^{(t)} \in \W_n \quad (t=0, \dots, t_n).
\]
Hence, for any $u>0$ we get
\begin{eqnarray*}
&&
\PROB \{ T_{2,n} > u \}
\\
&&
\leq
\PROB \Bigg\{
\exists f \in \F_n :
\EXP \left(
\left|
\frac{f(X)}{\beta_n} - \frac{T_{\beta_n}Y}{\beta_n}
\right|^2
\right)
-
\EXP \left(
\left|
\frac{m_{\beta_n}(X)}{\beta_n} - \frac{T_{\beta_n}Y}{\beta_n}
\right|^2
\right)
\\
&&\hspace*{3cm}-
\frac{1}{n} \sum_{i=1}^n
\left(
\left|
\frac{f(X_i)}{\beta_n} - \frac{T_{\beta_n}Y_i}{\beta_n}
\right|^2
-
\left|
\frac{m_{\beta_n}(X_i)}{\beta_n} - \frac{T_{\beta_n}Y_i}{\beta_n}
\right|^2
\right)
\Bigg\}
\\
&&\hspace*{2cm}
> \frac{1}{2} \cdot
\left(
\frac{u}{\beta_n^2}
+
\EXP \left(
\left|
\frac{f(X)}{\beta_n} - \frac{T_{\beta_n}Y}{\beta_n}
\right|^2
\right)
-
\EXP \left(
\left|
\frac{m_{\beta_n}(X)}{\beta_n} - \frac{T_{\beta_n}Y}{\beta_n}
\right|^2
\right)
\right),
\end{eqnarray*}
where
\[
\F_n = \left\{ T_{\beta_n} f_\bw \quad : \quad \bw \in \W_n \right\}.
\]
By Lemma \ref{le4} we get
\begin{eqnarray*}
&&
\Nu_1 \left(
\delta , \left\{
\frac{1}{\beta_n} \cdot f : f \in \F_n
\right\}
, x_1^n
\right)
\leq
\Nu_1 \left(
\delta \cdot \beta_n , \F_n
, x_1^n
\right)
\\
&&
\leq
\left(
\frac{ c_{24}}{\delta}
\right)^{
c_{25} \cdot (\log n)^{2d}  n^{\tau \cdot d} \cdot (\log n)^{2 \cdot (L-1) \cdot d}
\cdot
   \left(\frac{K_n \cdot (\log n)^2}{\beta_n \cdot \delta}\right)^{d/k} + c_{26}
  }.
\end{eqnarray*}
By choosing $k$ large enough we get for $\delta>1/n^2$
\[
\Nu_1 \left(
\delta , \left\{
\frac{1}{\beta_n} \cdot f : f \in \F_n
\right\}
, x_1^n
\right)
\leq
c_{27} \cdot n^{ c_{28} \cdot n^{\tau \cdot d + \epsilon/2}}.
\]
This together with Theorem 11.4 in Gy\"orfi et al. (2002) leads for $u
\geq 1/n$ to
\[
\PROB\{T_{2,n}>u\}
\leq
14 \cdot
c_{27} \cdot n^{ c_{28} \cdot n^{\tau \cdot d +  \epsilon/2}}
\cdot
\exp \left(
- \frac{n}{5136 \cdot \beta_n^2} \cdot u
\right).
\]
For $\epsilon_n \geq 1/n$ we can conclude
\begin{eqnarray*}
\EXP \{ T_{2,n} \}
& \leq &
\epsilon_n + \int_{\epsilon_n}^\infty \PROB\{ T_{2,n}>u \} \, du
\\
& \leq &
\epsilon_n
+
14 \cdot
c_{27} \cdot n^{ c_{28} \cdot n^{\tau \cdot d +
    \epsilon/2}}
\cdot
\exp \left(
- \frac{n}{5136 \cdot \beta_n^2} \cdot \epsilon_n
\right)
\cdot
\frac{5136 \cdot \beta_n^2}{n}.
\end{eqnarray*}
Setting
\[
\epsilon_n = \frac{5136 \cdot \beta_n^2}{n}
\cdot
c_{28}
\cdot
 n^{\tau \cdot d +
    \epsilon/2}
\cdot \log n
\]
yields the assertion of the fourth step of the proof.

In the {\it fifth step of the proof} we show
\begin{eqnarray*}
&&
\EXP \{ T_{4,n} \} \\
&&\leq
c_{29} \cdot \Bigg(
\sup_{
(\bar{w}_{i,j,k}^{(l)})_{i,j,k,l} :
\atop
|\bar{w}_{i,j,k}^{(l)}-w_{i,j,k}^{(l)}| \leq \log n
\, (l=0, \dots, L-1)
}
\int
|
\sum_{k=1}^{\tilde{K}_n}
w_{1,1,k} \cdot f_{\bar{\bw},k,j}^{(L)}(x)-m(x)|^2 \PROB_X (dx)
\\
&&
\hspace*{3.5cm}
+ \sum_{k=1}^{\tilde{K}_n} |w_{1,1,k}^{(L)}|^2
+
\frac{(\log n)^2}{n}
+
\frac{ n^{\tau \cdot d + \epsilon}}{n}
\Bigg).
\end{eqnarray*}
Using
\[
|T_{\beta_n} z - y| \leq |z-y|
\quad \mbox{for } |y| \leq \beta_n
\]
we get
\begin{eqnarray*}
&&
T_{4,n}/2
\\
&&
=
\Big[ \frac{1}{n} \sum_{i=1}^n
|m_n(X_i)-Y_i|^2
-
 \frac{1}{n} \sum_{i=1}^n
|m(X_i)- Y_i|^2
\Big] \cdot 1_{A_n}
\\
&&
\leq
\Big[
\frac{1}{n} \sum_{i=1}^n
|f_{\bw^{(t_n)}}(X_i)-Y_i|^2
-
 \frac{1}{n} \sum_{i=1}^n
|m(X_i)- Y_i|^2
\Big] \cdot 1_{A_n}
\\
&&
\leq
\big[ F_n(\bw^{(t_n)})
-
 \frac{1}{n} \sum_{i=1}^n
|m(X_i)- Y_i|^2
\Big] \cdot 1_{A_n}.
\end{eqnarray*}
Application of Lemma \ref{le1} (which is applicable on $A_n$ because
of Lemma \ref{le2} and Lemma \ref{le3}) implies that this in turn
is less than
\begin{eqnarray*}
&&
\big[
F_n(\bw^{(t_n-1)}) - \frac{1}{2 L_n} \cdot \| \nabla_\bw F_n(\bw^{(t_n-1)}) \|^2
-
 \frac{1}{n} \sum_{i=1}^n
|m(X_i)- Y_i|^2
\Big] \cdot 1_{A_n}.
\end{eqnarray*}
Since the sum of squares of all partial derivatives is at least
as large as the sum of squares of the partial derivatives with respect
to the outer weights $w_{1,1,k}^{(L)}$ $(k=1, \dots, K_n)$, we can
upper bound this in turn via Lemma \ref{le6} by
\begin{eqnarray*}
&&
\big[
F_n(\bw^{(t_n-1)}) - \frac{1}{2 L_n} \cdot 4 \cdot c_3 \cdot
  ( F_n(\bw^{(t_n-1)}) -  F_n((\bw^*)^{(t_n-1)})
\\
&&
\hspace*{6cm}
-
 \frac{1}{n} \sum_{i=1}^n
|m(X_i)- Y_i|^2
\Big] \cdot 1_{A_n}
\\
&&
=
\Big[
 \left(
1 - \frac{2 \cdot c_3}{ L_n}
\right)
\cdot
F_n(\bw^{(t_n-1)})
+
\frac{2 \cdot c_3}{ L_n} \cdot
F_n((\bw^*)^{(t_n-1)})
-
 \frac{1}{n} \sum_{i=1}^n
|m(X_i)- Y_i|^2
\Big] \cdot 1_{A_n}.
\end{eqnarray*}
Applying this argument repeatedly shows that
\begin{eqnarray*}
&&
T_{4,n}/2
\\
&&
\leq
\Big[
 \left(
1 - \frac{2 \cdot c_3}{ L_n}
\right)^{t_n}
\cdot
F_n(\bw^{(0)})
+
\sum_{k=1}^{t_n}
\frac{2 \cdot c_3}{ L_n} \cdot
  \left(
1 - \frac{2 \cdot c_3}{ L_n}
\right)^{k-1}
F_n((\bw^*)^{(t_n-k)})
\\
&&
\hspace*{3cm}
-
 \frac{1}{n} \sum_{i=1}^n
|m(X_i)- Y_i|^2
\Big] \cdot 1_{A_n}.
\end{eqnarray*}
This implies
\begin{eqnarray*}
&&
\EXP\{ T_{4,n}/2 \}
\\
&&
\leq
\left(
1 - \frac{2 \cdot c_3}{ L_n}
\right)^{t_n}
\cdot \EXP \{Y^2\}
+
\sum_{k=1}^{t_n}
\frac{2 \cdot c_3}{ L_n} \cdot
  \left(
1 - \frac{2 \cdot c_3}{ L_n}
\right)^{k-1}
\cdot
\\
&&
\hspace*{0.5cm}
\EXP \Bigg(
\Bigg(
\frac{1}{n} \sum_{i=1}^n
|
\sum_{k=1}^{{K}_n}
w_{1,1,k}^{(L)} \cdot f_{(\bw^*)^{(t_n-k)},j,1}^{(L)}(X_i)-Y_i|^2
-
\frac{1}{n} \sum_{i=1}^n
|m(X_i)-Y_i|^2
\Bigg) \cdot
   1_{A_n} \Bigg)
\\
&&
\quad
+ c_3 \cdot \sum_{k=1}^{\tilde{K}_n} |w_{1,1,k}^{(L)}|^2
\\
&&
\leq
\left(
1 - \frac{2 \cdot c_3}{ L_n}
\right)^{t_n}
\cdot \EXP \{Y^2\} + c_3 \cdot \sum_{k=1}^{\tilde{K}_n} |w_{1,1,k}^{(L)}|^2
\\
&&
\quad
+
\sum_{k=1}^{t_n}
\frac{2 \cdot c_3}{ L_n} \cdot
  \left(
1 - \frac{2 \cdot c_3}{ L_n}
\right)^{k-1}
\cdot 2 \cdot
\\
&&
\hspace*{1cm}
\Bigg(
\sup_{
(\bar{w}_{i,j,k}^{(l)})_{i,j,k,l} :
\atop
|\bar{w}_{i,j,k}^{(l)}-w_{i,j,k}^{(l)}| \leq \log n
\, (l=0, \dots, L-1)
}
\int
|
\sum_{k=1}^{\tilde{K}_n}
w_{1,1,k} \cdot f_{\bar{\bw},k,1}^{(L)}(x)-m(x)|^2 \PROB_X (dx)
\Bigg)
\\
&&
\quad
+
\sum_{k=1}^{t_n}
\frac{2 \cdot c_3}{ L_n} \cdot
  \left(
1 - \frac{2 \cdot c_3}{ L_n}
\right)^{k-1}
\cdot
\\
&&
\hspace*{0.5cm}
\EXP \Bigg(
\Bigg(
\frac{1}{n} \sum_{i=1}^n
|
\sum_{k=1}^{{K}_n}
w_{1,1,k} \cdot f_{(\bw^*)^{(t_n-k)},j,1}^{(L)}(X_i)-Y_i|^2
-
\frac{1}{n} \sum_{i=1}^n
|m(X_i)-Y_i|^2
\\
&&
\hspace*{0.5cm}
- 2 \cdot \left(
 \EXP\{ |\sum_{k=1}^{{K}_n}
w_{1,1,k} \cdot f_{(\bw^*)^{(t_n-k)},j,1}^{(L)}(X)-Y|^2 | \D_n\}
-\EXP\{|m(X)-Y|^2\}
\right)
\Bigg) \cdot
   1_{A_n} \Bigg),
\end{eqnarray*}
where the last inequality followed from (\ref{pth1eq2}).
Arguing as in the beginning of the proof (and using in particularly
the arguments from Steps 1, 2 and 4) we get
\begin{eqnarray*}
&&
\EXP \Bigg(
\Bigg(
\frac{1}{n} \sum_{i=1}^n
|
\sum_{k=1}^{{K}_n}
w_{1,1,k} \cdot f_{(\bw^*)^{(t_n-k)},j,1}^{(L)}(X_i)-Y_i|^2
-
\frac{1}{n} \sum_{i=1}^n
|m(X_i)-Y_i|^2
\\
&&
\hspace*{0.5cm}
- 2 \cdot \left(
 \EXP\{ |\sum_{k=1}^{{K}_n}
w_{1,1,k} \cdot f_{(\bw^*)^{(t_n-k)},j,1}^{(L)}(X)-Y|^2 | \D_n\}
-\EXP\{|m(X)-Y|^2\}
\right)
\Bigg) \cdot
   1_{A_n}
\Bigg)
\\
&&
\leq
c_{30} \cdot \frac{(\log n)^2}{n}
+
c_{31} \cdot
\frac{ n^{\tau \cdot d + \epsilon}}{n}.
\end{eqnarray*}

From this we
conclude
\begin{eqnarray*}
&&
\EXP\{ T_{4,n}/2 \} \\
&&
\leq
\left(
1 - \frac{2 \cdot c_3}{ L_n}
\right)^{t_n}
\cdot \EXP \{Y^2\}
\\
&&
\quad
+
4 \cdot
\Bigg(
\sup_{
(\bar{w}_{i,j,k}^{(l)})_{i,j,k,l} :
\atop
|\bar{w}_{i,j,k}^{(l)}-w_{i,j,k}^{(l)}| \leq \log n
\, (l=0, \dots, L-1)
}
\int
|
\sum_{k=1}^{\tilde{K}_n}
w_{1,1,k} \cdot f_{\bar{\bw},k,j}^{(L)}(x)-m(x)|^2
\PROB_X (dx)
\Bigg)
\\
&&
\quad
+ c_3 \cdot \sum_{k=1}^{\tilde{K}_n} |w_{1,1,k}^{(L)}|^2+
c_{30} \cdot \frac{(\log n)^2}{n}
+
c_{31} \cdot
\frac{ n^{\tau \cdot d + \epsilon}}{n}.
\end{eqnarray*}
The definition of $t_n$
together with (\ref{th1eq6}) implies
\[
\left(
1 - \frac{2 \cdot c_3}{ L_n}
\right)^{t_n}
\cdot \EXP \{Y^2\}
\leq
\exp \left(
-
\frac{2 \cdot c_3}{ L_n} \cdot t_n
\right)
\cdot \EXP \{Y^2\}
\leq \frac{c_{32}}{n}.
\]
Summarizing the above results we get the assertion.
    \hfill $\Box$

  \subsection{Auxiliary results for the proof of Theorem \ref{th2}}
    \label{se4sub3}

\begin{lemma}
  \label{le5}
Let $\sigma$ be the logistic squasher and let $0<\delta \leq 1$, $1 \leq \alpha_n \leq \log n$, $\mathbf{u},\mathbf{v} \in \mathbb{R}^d$ with
\begin{equation*}
v^{(l)}-u^{(l)} \geq 2\delta \quad \mbox{for } l\in\{1,\dots, d\}
\end{equation*}
and $x \in
  [-\alpha_n,\alpha_n]^d$.
Let $L, r,n, s \in \mathbb{N}$ with $L \geq 2$, $r \geq 2 \cdot d$,
$n \geq 8d$, $n\geq \exp(r+1)$ and $n \geq e^s$.
Let
\[
f_{\bw}(x)= f_{1,1}^{(L)}(x)
\]
where $f_{k,i}^{(l)}(x)$ are recursively defined by (\ref{se2eq2}) and (\ref{se2eq3}).

Assume
\begin{equation}
\label{le5eq1}
w_{1,j,j}^{(0)}=\frac{4d\cdot (\log n)^2}{\delta}
\quad \mbox{and} \quad  w_{1,j,0}^{(0)}=-\frac{4d \cdot (\log n)^2}{\delta} \cdot u^{(j)} \quad \text{for }j \in \{1, \dots, d\},
\end{equation}
\begin{equation}
\label{le5eq2}
w_{1,j+d,j}^{(0)}=-\frac{4d \cdot (\log n)^2}{\delta}
 \quad \mbox{and}\quad w_{1,j+d,0}^{(0)}=\frac{4d \cdot (\log n)^2}{\delta} \cdot v^{(j)}
\quad \text{for }j \in \{1, \dots, d\},
\end{equation}
\begin{equation}
\label{le5eq3}
w_{1,s,t}^{(0)}=0
\quad \mbox{if }
s \leq 2d,
s \neq t, s \neq t+d \mbox{ and } t>0,
\end{equation}
\begin{equation}
\label{le5eq4}
w_{1,1,t}^{(1)}= 8 \cdot (\log n)^2 \quad \mbox{for }t \in \{1, \dots, 2d\},
\end{equation}
\begin{equation}
\label{le5eq5}
w_{1,1,0}^{(1)} = - 8 (\log n)^2\left(2d-\frac{1}{2}\right),
\end{equation}
\begin{equation}
\label{le5eq6}
w_{1,1,t}^{(1)}=0
\quad \mbox{for }t > 2d,
\end{equation}
\begin{equation}
\label{le5eq7}
w_{1,1,1}^{(l)}= 6 \cdot (\log n)^2 \quad \mbox{for }l \in \{2, \dots, L\},
\end{equation}
\begin{equation}
\label{le5eq8}
w_{1,1,0}^{(l)} = - 3\cdot(\log n)^2 \quad \mbox{for }l \in \{2, \dots, L\}
\end{equation}
and
\begin{equation}
\label{le5eq9}
w_{1,1,t}^{(l)}=0
\quad \mbox{for }
 t > 1\mbox{ and } l \in \{2, \dots, L\}.
\end{equation}

Let $\bar{\bw}$ be such that
\begin{equation}
\label{le5eq10}
|\bar{w}_{1,i,j}^{(l)}-w_{1,i,j}^{(l)}| \leq \log n
\quad \mbox{for all } l=0,\dots,L-1.
\end{equation}
Then, we have
\[
f_{\bar{\bw}}(x) \geq 1 - \frac{1}{n^s} \mbox{ if } x \in
[u^{(1)}+\delta, v^{(1)}-\delta] \times \dots \times
[u^{(d)}+\delta, v^{(d)}-\delta]
\]
and
\[
f_{\bar{\bw}}(x) \leq \frac{1}{n^s} \mbox{ if } x^{(i)} \notin
[u^{(i)}-\delta, v^{(i)}+\delta] \mbox{ for some } i \in \{1,\dots, d\}.
\]
  \end{lemma}

\noindent
    {\bf Proof.}
The result follows from the proof of
Lemma 5 in Drews and Kohler (2022). In fact, in this proof it is shown that
\[
\sum_{j=1}^{r} \bar{w}_{1,1,j}^{(L-1)}\cdot \bar{f}_{1,j}^{(L-1)}(x) +
\bar{w}_{1,1,0}^{(l-1)}
\geq
2 (\log n)^2 - \frac{6}{n} (\log n)^2
\]
holds
if
$x \in
[u^{(1)}+\delta, v^{(1)}-\delta] \times \dots \times
[u^{(d)}+\delta, v^{(d)}-\delta]$,
and that
\[
\sum_{j=1}^{r} \bar{w}_{1,1,j}^{(L-1)}\cdot \bar{f}_{1,j}^{(L-1)}(x) +
\bar{w}_{1,1,0}^{(l-1)}
\leq
-2(\log n)^2 + 6 \log n \cdot \frac{1}{n}
\]
holds if $ x^{(i)} \notin
[u^{(i)}-\delta, v^{(i)}+\delta]$ for some $ i \in \{1,\dots, d\}$.
Because of $n \geq 6$ and $n \geq e^s$ we have
\[
2 (\log n)^2 - \frac{6}{n} (\log n)^2
\geq s \cdot \log n
\]
and
\[
-2(\log n)^2 + 6 \log n \cdot \frac{1}{n}
\leq - s \cdot \log n,
\]
from which we get the assertion as in the proof of
Lemma 5 in Drews and Kohler (2022).
    \hspace*{2cm} \hfill $\Box$

In our next lemma we use a multiscale approximation in order
to approximate a Lipschitz continuous function by a deep neural
network.

\begin{lemma}
\label{le8}
Let $1/2 \leq p \leq 1$, $C>0$,
let $f:\Rd \rightarrow \R$ be a $(p,C)$--smooth function and let
$X$ be an $\Rd$-valued random variable with $supp(X) \subseteq
[0,1]^d$.
Let $l \in \N$, $0<\delta<1/2$ with
\begin{equation}
\label{le8eq2}
c_{33} \cdot \delta \leq \frac{1}{2^l} \leq c_{34} \cdot \delta
\end{equation}
 and let $L,r,s \in \N$ with
\[
L \geq 2 \quad \mbox{and} \quad r \geq 2d
\]
and let
\[
\tilde{K}_n \geq \left( l \cdot (2^l+1)^{2d}+1 \right)^3
\]
Then there exist
\[
w_{k,i,j}^{(l)} \in
[-c_1 \cdot (\log n)^2, c_1 \cdot (\log n)^2]
\quad (l=1, \dots, L, k=1, \dots \tilde{K}_n)
\]
and
\[
w_{k,i,j}^{(0)}
\in \left[-\frac{8 \cdot d \cdot (\log n)^2}{\delta}, \frac{8 \cdot d \cdot
      (\log n)^2}{\delta} \right]
\quad (k=1, \dots, \tilde{K}_n).
\]
such that
for all $\bar{\bw}$ satisfying
$|\bar{w}_{i,j,k}^{(l)}-w_{i,j,k}^{(l)}| \leq \log n$
$(l=0, \dots, L-1)$ we have for $n$ sufficiently large
\begin{eqnarray}
\label{le8eq3}
&&
\int
|
\sum_{k=1}^{\tilde{K}_n}
w_{1,1,k}^{(L)} \cdot f_{\bar{\bw},k,1}^{(L)}(x)-f(x)|^2 \PROB_X (dx)
\nonumber \\
&&
\hspace*{2cm}
\leq
c_{35} \cdot \left(
l^2 \cdot \delta + \delta^{2p}
+
\frac{l \cdot (2^l+1)^{2d}}{n^s}
\right)
,
\end{eqnarray}
\begin{equation}
\label{le8eq1}
|
\sum_{k=1}^{\tilde{K}_n}
w_{1,1,k}^{(L)} \cdot f_{\bar{\bw},k,1}^{(L)}(x)| \leq
c_{36}
\cdot
\left( 1 +
\frac{(2^l+1)^{2d}}{n^s}
\right)
\quad (x \in [0,1]^d)
\end{equation}
and
\begin{equation}
\label{le8eq4}
\sum_{k=1}^{\tilde{K}_n} |w_{1,1,k}^{(L)}|^2 \leq \frac{c_{37}}{2^{2    \cdot d
    \cdot l}} .
\end{equation}
\end{lemma}

\noindent
{\bf Proof.}
In the proof we use a multiscale approximation of $f$ by piecewise
constant functions and use Lemma \ref{le5} in order to approximate
the piecewise constant functions by a linear combination of
neural networks.

In order to construct the multiscale approximation of $f$ by piecewise
constant functions we use a sequence of coverings
$\P^{(0)}=\{[0,1]^d\}$, $\P^{(1)}$, \dots, $\P^{(l)}$ of $[0,1]^d$
with the following properties:
\begin{enumerate}
\item
$\P^{(k)}$ consists of $(2^k+1)^d$ many pairwise disjoint cubes of side length $1/2^k$
$(k=1, \dots, l)$.

\item
$[0,1]^d \subseteq \cup_{A \in \P^{(k)}} A$

\item
\begin{equation}
\label{ple8eq1}
\PROB_X \left(
\cup_{A \in \P^{(k)}} A_{border,\delta}
\right)
\leq
4 d \cdot 2^k \cdot \delta,
\end{equation}
where
\begin{eqnarray*}
A_{border,\delta}
&=&
[u^{(1)}-\delta, v^{(1)}+\delta] \times \dots \times
[u^{(d)}-\delta, v^{(d)}+\delta]
\\
&&
\hspace*{3cm}
\setminus
[u^{(1)}+\delta, v^{(1)}-\delta] \times \dots \times
[u^{(d)}+\delta, v^{(d)}-\delta]
\end{eqnarray*}
for
\[
A=[u^{(1)}, v^{(1)}] \times \dots \times
[u^{(d)}, v^{(d)}].
\]
\end{enumerate}

We can ensure (\ref{ple8eq1}) by shifting a partition of
\[
\left[ - \frac{1}{2^k}, 1 \right]^d
\]
consisting of $(2^k+1)^d$ many cubes of side length $1/2^k$
separately
in each component by multiples of $2 \cdot \delta$ less than or
equal to $1/2^k$, which gives us for each component
\[
\left\lfloor \frac{1}{2 \cdot \delta} \cdot \frac{1}{2^k} \right \rfloor
\]
disjoint sets of which at least one must have $\PROB_X$-measure
less than or equal to
\[
\frac{1}{\lfloor \frac{1}{2 \cdot\delta} \cdot \frac{1}{2^k} \rfloor}
\leq
\frac{1}{ \frac{1}{2 \cdot\delta} \cdot \frac{1}{2^k} -1 }
\leq
\frac{2 \cdot \delta \cdot 2^k}{1-2 \cdot \delta \cdot 2^k} \leq
4 \cdot \delta \cdot 2^k
\]
in case $2 \cdot \delta \cdot 2^k \leq 1/2$, which we can assume
w.l.o.g.

For $x \in [0,1]^d$ denote by $z_{\P^{(k)}}(x)$ the center $z_A$ of
the
set $A \in \P^{(k)}$ which contains $x$.

We approximate
\[
f(x)
\]
by
\begin{eqnarray}
\label{ple8eq3}
&&f(z_{\P^{(l)}}(x))
\nonumber \\
&&=
f(z_{\P^{(0)}}(x))+
\sum_{k=1}^l \left(
f(z_{\P^{(k)}}(x))-f(z_{\P^{(k-1)}}(x))
\right)
\nonumber \\
&&=
f(z_{\P^{(0)}}(x))+
\sum_{k=1}^l
\sum_{A_1 \in \P^{(k)}, A_2 \in \P^{(k-1)}: A_1 \cap A_2 \neq
    \emptyset}
\left(
f(z_{A_1}) - f(z_{A_2})
\right) \cdot 1_{A_1 \cap A_2}(x).
\end{eqnarray}

For a $d$-dimensional rectangle $R$ let
\[
f_{net,R,\delta}
\]
be the neural network from Lemma \ref{le5} which approximates
\[
1_R.
\]
(In case that $2\delta$ is less than the minimal side length of $R$,
Lemma \ref{le5} does not imply that $f_{net,R,\delta}$ is in the
inner part of $R$ close to one).

We approximate $f(z_{\P^{(l)}}(x))$ by
\begin{eqnarray*}
&&
f(z_{\P^{(0)}}(x)) \cdot f_{net, [-1,2]^d, \delta}(x)
\\
&&
\quad
+
\sum_{k=1}^l
\sum_{A_1 \in \P^{(k)}, A_2 \in \P^{(k-1)}: A_1 \cap A_2 \neq
    \emptyset}
\left(
f(z_{A_1}) - f(z_{A_2})
\right) \cdot f_{net,A_1 \cap A_2,\delta}(x)
\end{eqnarray*}
and let $\bw$ be the weights of the above neural network.
Observe that this network consists of
\[
1 +
\sum_{k=1}^l |\{A_1 \in \P^{(k)}, A_2 \in \P^{(k-1)} \, : \,A_1 \cap A_2 \neq
    \emptyset\}|
\leq
1 + l \cdot (2^l +1)^{2d}
\]
many fully connected neural networks which are computed in parallel.
Denote the neural network where the weights
of $f_{net,R,\delta}$ are replaced by the corresponding weights of
$\bar{\bw}$
by
$f_{net,\bar{\bw},R,\delta}$. We then set
\begin{eqnarray*}
f_{net}(x)
&=&
\sum_{k=1}^{\tilde{K}_n}
w_{1,1,k}^{(L)} \cdot f_{\bar{\bw},k,1}^{(L)}(x)
\\
&
=&
f(z_{\P^{(0)}}(x)) \cdot f_{net, \bar{\bw},[-1,2]^d, \delta}(x)
\\
&&
+
\sum_{k=1}^l
\sum_{A_1 \in \P^{(k)}, A_2 \in \P^{(k-1)}: A_1 \cap A_2 \neq
    \emptyset}
\left(
f(z_{A_1}) - f(z_{A_2})
\right) \cdot f_{net,\bar{\bw},A_1 \cap A_2,\delta}(x)
.
\end{eqnarray*}
Since we have $\delta<1/2$ we know by Lemma \ref{le5} that
$f_{net, \bar{\bw},[-1,2]^d, \delta}(x) \geq 1-1/n^s$ holds for all $x \in [0,1]^d$.
We have
\begin{eqnarray*}
&&
\int
| f_{net}(x)-f(x)|^2 \PROB_X (dx)
\\
&&
\leq
2 \cdot
\int
| f(z_{\P^{(l)}}(x))-f(x)|^2 \PROB_X (dx)
+
2 \cdot
\int
| f_{net}(x)-f(z_{\P^{(l)}}(x))|^2 \PROB_X (dx)
.
\end{eqnarray*}
By Lemma \ref{le5}, (\ref{ple8eq3}) and
the $(p,C)$--smoothness of $f$, which implies
\[
| f(z_{\P^{(l)}}(x)) - f(x)|
\leq  c_{38} \cdot \left( \frac{1}{2^l} \right)^p
\quad \mbox{for all }
x \in [0,1]^d
\]
and
\begin{equation}
\label{ple8eq2}
|f(z_{A_1}) - f(z_{A_2})| \leq  c_{39} \cdot \left( \frac{1}{2^k} \right)^p
\quad \mbox{for all }
A_1 \in \P^{(k)}, A_2 \in \P^{(k-1)} \mbox{ with } A_1 \cap A_2 \neq \emptyset,
\end{equation}
we can bound the last sum above by
\begin{eqnarray*}
&&
c_{40} \cdot \left( \frac{1}{2^l} \right)^{2p}
+ c_{41} \cdot \left( \frac{1}{n^s} \right)^2
\\
&&
\hspace*{1cm}
+
4 \cdot l \cdot \sum_{k=1}^l
\int
\big|
\sum_{A_1 \in \P^{(k)}, A_2 \in \P^{(k-1)}: A_1 \cap A_2 \neq
    \emptyset}
\left(
f(z_{A_1}) - f(z_{A_2})
\right) \\
&&
\hspace*{4cm}
\cdot (1_{A_1 \cap A_2}(x) - f_{net,A_1 \cap A_2,\delta}(x))
\big|^2 \PROB_X (dx)
\\
&&
\leq
c_{40} \cdot \left( \frac{1}{2^l} \right)^{2p}
+ c_{41} \cdot \left( \frac{1}{n^s} \right)^2
+
4 \cdot l \cdot \sum_{k=1}^l
\Big(
c_{42} \cdot
\left( \frac{1}{2^k} \right)^{2p}
\cdot \PROB_X( \cup_{A \in \P^{(k)} \cup \P^{(k-1)}} A_{border, \delta})
\\
&&
\hspace*{1cm}
+ c_{43} \cdot \left( \left(\frac{1}{2^k} \right)^p \cdot \frac{1}{n^s} \right)^2 \cdot (2^k+1)^{2d}
\Big)
\\
&&
\leq
c_{44} \cdot \left(
\left( \frac{1}{2^l} \right)^{2p} +
l^2 \cdot \delta + l \cdot (2^l+1)^{2d} \cdot \frac{1}{n^s}
\right).
\end{eqnarray*}
By (\ref{le8eq2}) this implies (\ref{le8eq3}).

Next we prove (\ref{le8eq1}).
By construction we know that $f_{net,R,\delta}$ is bounded
in absolute value by one, which implies
\begin{eqnarray*}
  &&
  \left|
\sum_{k=1}^{\tilde{K}_n}
w_{1,1,k}^{(L)} \cdot f_{\bar{\bw},k,1}^{(L)}(x)
\right|
\\
&
&\leq
|f(z_{\P^{(0)}}(x))|
+
\sum_{k=1}^l
\sum_{A_1 \in \P^{(k)}, A_2 \in \P^{(k-1)}: A_1 \cap A_2 \neq
    \emptyset}
\left|
f(z_{A_1}) - f(z_{A_2})
\right| \cdot
f_{net,\bar{\bw}, A_1 \cap A_2, \delta}(x).
\end{eqnarray*}
Using
(\ref{ple8eq2}),
Lemma \ref{le6} and that for each
$k \in \{1, \dots, l\}$ there are at most $c_{45}$
many $A_1 \in \P^{(k)}$, $A_2 \in \P^{(k-1)}$
such that
\[
A_1 \cap A_2 \neq \emptyset
\quad \mbox{and} \quad
x \in (A_1 \cap A_2) \cup
(A_1 \cap A_2)_{border, \delta},
\]
we can bound the term on the right-hand side above by
\begin{eqnarray*}
  &&
\|f\|_\infty + \sum_{k=1}^l \left(
c_{45} \cdot c_{39} \cdot
\left( \frac{1}{2^k} \right)^p
+
c_{39} \cdot
\left( \frac{1}{2^k} \right)^p
\cdot
(2^k+1)^{2d} \cdot \frac{1}{n^s}
\right)
\\
&&
\leq
c_{46} \cdot
\left(
1
+
\frac{
(2^l+1)^{2d}
}{n^s}
\right).
  \end{eqnarray*}

There at most $ l \cdot (2^l+1)^{2d}+1$ many output weights of the
above neural network are all bounded in absolute value by a constant,
which implies that if we use $\tilde{K}_n =l \cdot (2^l+1)^{2d}+1$
we will get
\[
\sum_{k=1}^{\tilde{K}_n} |w_{1,1,k}^{(L)}|^2 \leq c_{47} \cdot \tilde{K}_n.
\]
In order to get a smaller upper bound, we repeat the whole
construction $( l \cdot (2^l+1)^{2d}+1)^2$ many times, each time with output
weights divided by $ (l \cdot (2^l+1)^{2d}+1)^2$. The above proof
implies
that the linear combination of these $ (l \cdot (2^l+1)^{2d}+1)^3$ many
neural networks still satisfies (\ref{le8eq3}) and (\ref{le8eq1}) (here
we use that in each of the networks we can use the same exception sets
where $f_{net,R,\delta}$ is not accurate). This results in
\[
\sum_{k=1}^{\tilde{K}_n} |w_{1,1,k}^{(L)}|^2 \leq
\sum_{k=1}^{(l \cdot (2^l+1)^{2d}+1)^3 }
\left(\frac{c_{47}}{  (l \cdot (2^l+1)^{2d}+1)^2}\right)^2
\leq
\frac{c_{48}}{2^{2d \cdot l}}.
\]
\quad \hfill $\Box$

\subsection{Proof of Theorem \ref{th2}}
\label{se4sub4}
Set
\[
l=\lfloor \frac{1}{1+d} \cdot \log n \rfloor
\]
(which implies that
that
(\ref{le8eq2}) holds for $\delta= n^{-1/(1+d)}$) and $\tilde{K}_n= (l \cdot
(2^l+1)^{2d}+1)^3$.
By applying Lemma \ref{le8} with a sufficiently large $s$
together with Theorem \ref{th1}
we get for $n$ sufficiently large
\begin{eqnarray*}
&&
\EXP \int | m_n(x)-m(x)|^2 \PROB_X (dx)
\leq
c_7 \cdot
\Bigg(
\frac{ n^{\frac{1}{1+d} \cdot d + \epsilon}}{n} +\sum_{k=1}^{\tilde{K}_n} |w_{1,1,k}^{(L)}|^2
\\
&&
\hspace*{1cm}
+
\sup_{
(\bar{w}_{i,j,k}^{(l)})_{i,j,k,l} :
\atop
|\bar{w}_{i,j,k}^{(l)}-w_{i,j,k}^{(l)}| \leq \log n
\, (l=0, \dots, L-1)
}
\int
|
\sum_{k=1}^{\tilde{K}_n}
w_{1,1,k}^{(L)} \cdot f_{\bar{\bw},k,1}^{(L)}(x)-m(x)|^2 \PROB_X (dx)
\Bigg)
\\
&&
\leq
c_{49} \cdot
\left(
\frac{ n^{\frac{1}{1+d} \cdot d + \epsilon}}{n} + n^{-\frac{2d}{d+1}} + (\log
   n)^2 \cdot  n^{-1/(1+d)}
+
n^{-\frac{2p}{d+1}}
\right)
\\
&&
\leq
c_{50} \cdot
n^{
- \frac{1}{1+d} + \epsilon
}.
\end{eqnarray*}

\hfill $\Box$

\subsection{Auxiliary results for the proof of Theorem \ref{th3}}
\label{se4sub5}

         \begin{lemma}
          \label{le7}
 Let $\alpha \geq 1$, $\beta>0$ and let $A,B,C \geq 1$.
  Let $\sigma:\R \rightarrow \R$ be $k$-times differentiable
  such that all derivatives up to order $k$ are bounded on $\R$.
  Let $\F$
  be the set of all functions
\[
f_{\bw}(x) = \sum_{I \subseteq \{1, \dots, d\} \, : \, |I|=d^* } f_{\bw_I} (x_I)
\]
where $f_{\bw_I}$ as defined by (\ref{se2eq1})--(\ref{se2eq3}) with
$d$ replaced by $d^*$ and weight vector $\bw_I$,
\[
\bw= \left( \bw_I \right)_{I \subseteq \{1, \dots, d\} \, : \, |I|=d^* },
\]
and where for any
$I \subseteq \{1, \dots, d\}$ with $|I|=d^*$
the weight vector $\bw_I$ satisfies
  \begin{equation}
    \label{le7eq1}
    \sum_{j=1}^{K_n} |(\bw_I)_{1,1,j}^{(L)}| \leq C,
    \end{equation}
  \begin{equation}
    \label{le7eq2}
    |(\bw_I)_{k,i,j}^{(l)}| \leq B \quad (k \in \{1, \dots, K_n\},
    i,j \in \{1, \dots, r\}, l \in \{1, \dots, L-1\})
    \end{equation}
and
  \begin{equation}
    \label{le7eq3}
    |(\bw_I)_{k,i,j}^{(0)}| \leq A \quad (k \in \{1, \dots, K_n\},
    i \in \{1, \dots, r\}, j \in \{1, \dots,d\}).
  \end{equation}
  Then we have for any $1 \leq p < \infty$, $0 < \epsilon < \beta$ and
  $x_1^n \in [-\alpha,\alpha]^d$
  \begin{eqnarray*}
    &&
  \Nu_p \left(
\epsilon, \{ T_\beta f  \, : \, f \in \F \}, x_1^n
\right)
\\
&&
\leq
\left(
c_{51} \cdot \frac{\beta^p }{\epsilon^p}
\right)^{
c_{52} \cdot \alpha^{d^*} \cdot A^{d^*} \cdot B^{(L-1) \cdot d^*} \left(\frac{C}{\epsilon}\right)^{d^*/k} + c_{53}
  }.\\
  \end{eqnarray*}
       \end{lemma}

         \noindent
             {\bf Proof.}
By the proof of Lemma \ref{le4} we have for any
$I \subseteq \{1, \dots, d\}$ with $|I|=d^*$,
$x \in \R^{d^*}$ and any $s_1, \dots, s_k \in \{1, \dots, d\}$
\[
\left|
\frac{
\partial^k f_{\bw_I}
}{
\partial x^{(s_1)} \dots \partial x^{(s_k)}
}
(x)
\right|
\leq
c_{54} \cdot C \cdot B^{(L-1) \cdot k} \cdot A^k=:c.
\]
For
$I \subseteq \{1, \dots, d\}$ with $|I|=d^*$
let $\G \circ \Pi_I$ be the set of all piecewise polynomials
of total degree less than $k$ with respect to a partition $\Pi_I$
of $[-\alpha,\alpha]^{d^*}$ into cubes of sidelength
\[
\left(
c_{55} \cdot
\frac{\epsilon}{c}
\right)^{1/k},
\]
where $c_{55}=c_{55}(d,k)$ is a suitable constant greater zero.
Then a standard bound on the remainder of a multivariate Taylor
polynomial
shows that for each $f_{\bw_I}$ we can find $g_{\bw_I} \in \G \circ \Pi_I$ such
that
\[
|f_{\bw_I}(x) - g_{\bw_I}(x)| \leq \frac{\epsilon}{2 \cdot \left( d \atop d^* \right)}
\]
holds for all $x \in [-\alpha,\alpha]^{d^*}$.
Let $\HH$ be the set of all functions of the form
\[
h(x) = \sum_{I \subseteq \{1, \dots, d\} \, : \, |I|=d^* } g_{\bw_I} (x_I)
\quad (x \in [-\alpha,\alpha]^d)
\]
$(g_I \in \G \circ \Pi_I, I \subseteq \{1, \dots, d\}, |I|=d^* )$.
Then we have
\[
  \Nu_p \left(
\epsilon, \{ T_\beta f  \, : \, f \in \F \}, x_1^n
\right)
\leq
  \Nu_p \left(
\frac{\epsilon}{2}, \{ T_\beta h  \, : \, h \in \HH \}, x_1^n
\right).
\]
$\HH$ is a linear vector space of dimension less than or equal to
\[
c_{56} \cdot
\left( d \atop d^* \right)
\cdot
\alpha^{d^*} \cdot \left(
\frac{c}{\epsilon}
\right)^{d^*/k},
\]
from which we get the assertion by an application of Theorems 9.4 and 9.5
in Gy\"orfi et al. (2002).
    \hfill $\Box$

\subsection{Proof of Theorem \ref{th3}}
\label{se4sub6}

Set
\[
l=\lfloor \frac{1}{1+d^*} \log n \rfloor
\]
(which implies that
that
(\ref{le8eq3}) holds for $\delta= n^{-1/(1+d^*)}$  and $\tilde{K}_n= (l \cdot
(2^l+1)^{2d}+1)^3$.

By applying Lemma \ref{le8} with a sufficiently large $s$
together with Theorem \ref{th1} and Lemma \ref{le7}
we get
\begin{eqnarray*}
&&
\EXP \int | m_n(x)-m(x)|^2 \PROB_X (dx)
\leq
c_7 \cdot
\Bigg(
\frac{ n^{\frac{1}{1+d} \cdot d + \epsilon}}{n} +\sum_{k=1}^{\tilde{K}_n} |w_{1,1,k}^{(L)}|^2
\\
&&
\hspace*{1cm}
+
\sup_{
(\bar{w}_{i,j,k}^{(l)})_{i,j,k,l} :
\atop
|\bar{w}_{i,j,k}^{(l)}-w_{i,j,k}^{(l)}| \leq \log n
\, (l=0, \dots, L-1)
}
\int
|
\sum_{k=1}^{\tilde{K}_n}
w_{1,1,k}^{(L)} \cdot f_{\bar{\bw},k,j}^{(L)}(x)-m(x)|^2 \PROB_X (dx)
\Bigg)
\\
&&
\leq
c_{57} \cdot
\left(
\frac{ n^{\frac{1}{1+d} \cdot d + \epsilon}}{n} + n^{- \frac{2p}{d^*+1}} + (\log
   n)^2 \cdot  n^{-1/(1+d)}
\right)
\\
&&
\leq
c_{58} \cdot
n^{
- \frac{1}{1+d} + \epsilon
}.
\end{eqnarray*}

\hfill $\Box$

\end{document}